\documentstyle[amscd,amssymb,12pt]{amsart}

\newtheorem{Thm}{Theorem}[section]

\newtheorem{Cor}[Thm]{Corollary}
\newtheorem{Lem}[Thm]{Lemma}
\newtheorem{Prop}[Thm]{Proposition}
\newtheorem{Def}[Thm]{Definition}
\newtheorem{Rmk}[Thm]{Remark}

\newcommand{\e}{\varepsilon}

\newcommand{\supp}{{\rm supp}\,}
\newcommand{\spn}{{\rm span}\,}

\def \N{\rm {\bf N}}
\def \R{\rm {\bf R}}

\begin{document}

\title[Weak Hilbert]%
      {Some more  weak Hilbert spaces}

\author{George Androulakis}
\address{Department of Mathematics, Math. Sci. Bldg.,
University of Missouri-Columbia, Columbia MO  65211}
\email{giorgis@@math.missouri.edu}

\author{Peter G. Casazza}
\address{Department of Mathematics, Math. Sci. Bldg.,
University of Missouri-Columbia, Columbia MO  65211}
\email{pete@@casazza.math.missouri.edu}
\thanks{The second author was supported by NSF DMS 9706108.}

\author{Denka N. Kutzarova}
\address{Current:  Dept. Math. \& Stat., Miami University,
Oxford, Ohio,  Permanent:  Institute of Mathematics,
Bulgarian Academy of Sciences, 1113 Sofia, Bulgaria}
\email{denka@@math.acad.bg}
\thanks{The third author  was partially supported by
the Bulgarian Ministry of Education and Science under contract MM-703/97}

\date{}
\maketitle

\bigskip
\bigskip
\noindent
{\bf Abstract:}
We give new examples of weak Hilbert spaces.

\section{Introduction}
The Banach space properties {\bf weak type 2} and {\bf weak cotype  
2} were
introduced and studied by V. Milman and G. Pisier \cite{MP}.   
Later, Pisier
\cite{P1} studied spaces which are both of weak type 2 and weak  
cotype 2 and called them {\bf weak Hilbert spaces}.  Weak Hilbert  
spaces are stable under passing to subspaces, dual spaces, and   
quotient spaces.  The canonical
example of a weak Hilbert space which is not a Hilbert space is  
convexified
Tsirelson space $T^2$ \cite{CS,J1,J2,P1}.  Tsirelson's space was
introduced by B.S. Tsirelson \cite{T} as the first example of a  
Banach space
which does not contain an isomorphic copy of $c_{0}$ or  
${\ell}_{p}$, $1\le p<
\infty$.  Today, we denote by $T$ the dual space of the original  
example of
Tsirelson since in $T$ we have an important analytic description of the
norm due to Figiel and Johnson \cite{FJ}.  In \cite{J1}, Johnson  
introduced
{\bf modified Tsirelson space} $T_{M}$.  Later, Casazza and Odell  
\cite{CO}
proved the surprising fact that $T_{M}$ is naturally isomorphic to the
original Tsirelson space $T$.  At this point, all the non-trivial  
examples of weak Hilbert spaces (i.e. those which are not Hilbert  
spaces) had unconditional
bases and had subspaces which failed to contain ${\ell}_{2}$.  A.  
Edgington
\cite{E} introduced a class of weak Hilbert spaces with  
unconditional bases
which are ${\ell}_{2}$-saturated.  That is, every subspace of the space
contains a further subspace isomorphic to a Hilbert space but the space
itself is not isomorphic to a Hilbert space.  R. Komorowski  
\cite{K} (or more
generally Komorowski and Tomczak-Jaegermann \cite{KT}) proved that
there are weak Hilbert spaces with no unconditional basis.  In fact, they
show that $T^2$ has such subspaces.  In another surprise, Nielsen and
Tomczak-Jaegermann \cite{NTJ} showed that all weak Hilbert spaces with
unconditional bases are very much like $T^{(2)}$.

There are still many open questions concerning weak Hilbert spaces and 
$T^{(2)}$, due partly to the shortage of non-trivial examples in this area.
For example, it is still a major open question in
the field whether a Banach space for which every subspace has an  
unconditional
basis (or just local unconditional structure - LUST) must be  
isomorphic to a Hilbert
space.  If there are such examples, they will probably come from the
class of weak Hilbert spaces. It is an open question whether every
weak Hilbert space has a basis, although Maurey and Pisier (see
\cite{M}) showed that separable weak Hilbert spaces have finite
dimensional decompositions. Nielsen and Tomczak-Jaegermann have shown
that weak Hilbert spaces that are Banach lattices have the property
that every subspace of every quotient space has a basis. But it
is unknown whether every weak Hilbert space can be  
embedded into
a weak Hilbert space with a unconditional basis.  In fact, it is unknown
if a weak Hilbert space embeds into a Banach lattice of finite  
cotype.  It turns out that this question is equivalent to the  
question of whether
every subspace of a weak Hilbert space must have the GL-Property  
\cite{CN}
which is slightly weaker than having LUST.  In
this note we extend the list of non-trivial examples of weak  
Hilbert spaces
by producing examples which are ${\ell}_{2}$-saturated but have no  
subspaces
isomorphic to subspaces of the previously known examples.

\section{Basic Constructions}
If $F$ is a finite dimensional Banach space then let $d(F)$ denote the
Banach-Mazur distance between $F$ and $\ell_2^{{\rm dim}F}$. The
fundamental notion of this note is the one of the weak Hilbert
space. Recall the following definition as one of the many equivalent
ones (cf \cite{P1} Theorem 2.1).

\begin{Def} \label{weakH}
A Banach space $X$ is said to be a weak Hilbert space if there exist
$\delta >0$ and $C \geq 1$ such that for every finite dimensional
subspace $E$ of $X$ there exists a subspace $F \subseteq E$ and a
projection $P:X \to F$ such that ${\rm dim} F \geq \delta {\rm dim}
E$, $ d(F) \leq C$ and $\| P \| \leq C$.
\end{Def}

\noindent
We need to recall the definition of the Schreier sets $S_n$, $n \in \N$
\cite{AA}. For $F, G \subset \N$,  we write $F < G$ when
$\max (F) < \min (G)$ or one of them is empty, and we write $n \leq F$
instead
of $\{ n \} \leq F$. 
$$
S_0 = \{ \{ n \} : n \in \N \} \cup \{ \emptyset \}.
$$ 
If $n \in \N \cup \{ 0 \}$ and $S_n$ has been defined, 
$$
S_{n+1} = \{ \cup_1 ^n F_i : n \in \N , 
n \leq F_1 <F_2< \cdots <F_n \mbox{ and } F_i \in S_n \mbox{ for } 
1 \leq i \leq n \}. 
$$ 
For $n \in \N$ a family of finite non-empty subsets $(E_i)$
of $\N$ is said to be {\em $S_n$-admissible} if $E_1<E_2< \cdots $ and
$(\min(E_i))\in S_n$. Also, $(E_i)$  is said to be
{\em $S_n$-allowable} if $E_i \cap E_j = \emptyset$ for $i \not = j$
and $(\min(E_i))\in S_k$.

\bigskip
\noindent
Every Banach space with a basis can be viewed as the completion of $c_{00}$
(the linear space of finitely supported real valued sequences) 
under a certain norm. $(e_i)$
will denote the unit vector basis for $c_{00}$ and whenever a Banach space 
$(X, \| \cdot \|)$ with
a basis is regarded as the completion of $(c_{00}, \| \cdot \|)$,
  $(e_i)$ will denote this (normalized)  basis.
If $x \in c_{00}$ and $E \subseteq \N$, $Ex \in c_{00}$ is the restriction 
of $x$ to $E$;
$(Ex)_j=x_j$ if $j \in E$ and 0 otherwise.
Also the {\em support} \/ of $x$, $\supp(x)$, (w.r.t. $(e_i)$) is the set 
$\{ j \in \N : x_j \not = 0\}$. If $f:\R \to \R$ is a function with $f(0)=0$ then  
for $x \in c_{00}$ $f(x)$ will denote the vector $f(x)= (f(x_i))$ in $c_{00}$.

\bigskip
\noindent
Let $(X, \| . \|)$ be a Banach space with an unconditional basis.
The norm of $X$ is $2$-convex provided that 
$$ 
\| (x^2 + y^2)^{1/2} \| \leq ( \| x \| ^2 + \| y \|^2 )^{1/2}
$$ 
for all vectors $x, y \in X$. The $2$-convexification of $(X, \| . \|)$ is
the Banach space $(X^{(2)}, \| . \|_{(2)})$ with an unconditional
basis, where $x \in X^{(2)}$ if and only if $x^2 \in X$ and
$$
\| x \|_{(2)} = \| x^2 \|^{1/2}.
$$   
Of course $\| . \|_{(2)}$ is $2$-convex. For $C>0$, $1 \leq p <
\infty$,  $n \in \N$ and $x_1, x_2,
\ldots , x_n \in X$ we say that $(x_i)_{i=1}^n$ is $C$-equivalent to
the unit vector basis of $\ell_p^n$ if there exist constants $A,B >0$
with $AB \leq C$ such that 
$$
\frac{1}{A}  (\sum |a_i|^p)^{1/p} \leq \| \sum a_i x_i \| \leq 
B (\sum |a_i|^p)^{1/p}
$$
for every sequence of scalars $(a_i)_{i=1}^n$. 
For $C>0$, we say that $X$ is  an {\em asymptotic $\ell_p$ space} 
(resp. {\em asymptotic $\ell_p$ space for vectors with
disjoint supports}) {\em with constant $C$} if  for every
$n$ and for every sequence of vectors $(x_i)_{i=1}^n$ such that  $(\supp
(x_i))_{i=1}^n$ is $S_1$-admissible (resp. $S_1$-allowable), we have
that  $(x_i)_{i=1}^n$ are
$C$-equivalent to the unit vector basis of $\ell_p^n$.

\bigskip
\noindent
If $(\| . \|)_n$ is a sequence of norms in $c_{00}$ then 
$\Sigma(\|. \|_n)$ will denote the completion of $c_{00}$  under the
norm 
$$
\| x \|_{\Sigma(\| . \|_n)}= \sum_{n=1}^\infty \| x \|_n.
$$

\bigskip
\noindent
Fix a sequence $\alpha= (\alpha_n)_{n \in \N}$ of elements of  $(0,1)$ and
$\ell, u \in (0,1)$ with $0< \ell \leq \frac{\alpha_{n+1}}{\alpha_n}
\leq u <1$  for all $n$ and $\sum_n \alpha_n =1$ (the existence of
numbers  $\ell, u \in (0,1)$ such that
the last relationships are valid will always be assumed whenever a
sequence $(\alpha_n)$ will be considered in these notes).
Edgington defined a sequence of norms $(\|. \|_{E,n})$
on $c_{00}$ by
$$
\| x \|_{E,0}= \| x \|_\infty, \quad \| x \|_{E, n+1}^2 = \sup \{ \sum_i \|
E_i x \|_{E,n}^2 : (E_i)_i \mbox{ is }S_1\mbox{-admissible} \}.
$$
Then Edgington defined the norm $\| . \|_{E_\alpha}$ by 
$$
\| x \|_{E_\alpha}= \left( \sum_n \alpha_n \| x \|_{E,n}^2 \right)^{1/2}.
$$
Let $E_\alpha$ denote the completion of $c_{00}$ with respect to $\|
. \|_E$. 
It is shown in \cite{E} that $E_\alpha$ is a weak Hilbert space which
is not  isomorphic to $\ell_2$, yet it is $\ell_2$-saturated. It is
easy to see that the spaces constructed by Edgington are asymptotic
$\ell_2$ spaces for vectors with disjoint supports. The main theorem
that we prove in these notes  (Theorem \ref{main}) shows that such
spaces are weak Hilbert spaces. 

\bigskip
\noindent
Let $(|.|_n)_{n \in \N}$ denote the sequence of the Schreier norms on $c_{00}$:
$$
| x |_n = \sup_{S \in S_n} \sum_{j \in S} |x_j|
$$
(if $x = \sum_j x_j e_j$). Then the weak Hilbert space $E_\alpha$ that was
constructed by Edgington \cite{E} is the $2$-convexification of $\Sigma
(\alpha_n |.|_n)$. One can see that $\Sigma( \alpha_n |.|_n)$ is an
asymptotic $\ell_1$ space for vectors with disjoint supports which is
$\ell_1$-saturated, yet not isomorphic to $\ell_1$. In these notes we
give examples of sequences of norms that can replace $(|.|_n)$ in
$\Sigma ( \alpha_n |.|_n)$ to obtain asymptotic $\ell_1$ spaces for
vectors with disjoint supports which are $\ell_1$-saturated yet not
isomorphic to $\ell_1$. The $2$-convexification of each of these
spaces will give $\ell_2$ saturated weak Hilbert spaces which are not
isomorphic to $\ell_2$.

\bigskip
\noindent
{\em Definition of the spaces $V$, $W$, $V'$ and $W'$:}
Let $(\theta_n)_{n \in \N}$ be a sequence of real
numbers in $(0,1)$ with $\lim_n \theta_n=0$ (this assumption will always be valid
whenever a sequence $(\theta_n)$ will be considered in these notes) and
let $s \in \N$. The asymptotic $\ell_1$ spaces  
$V=T_M(\theta_n, S_n)_n$ and $W=T_{M(s)}(\theta_n,S_n)_n$ were introduced
in \cite{AD} and \cite{ADKM} as the completion of $c_{00}$ under the
norms:
\begin{eqnarray*}
\| x \|_V & = & \| x \|_\infty \vee \sup_n 
  \sup \{ \theta_n \sum_i \| E_i x \|_V : (E_i) \mbox{ is }S_n 
\mbox{ allowable} \},\\
\| x \|_W & = & \| x \|_\infty \vee \sup_{n \leq s} 
  \sup \{ \theta_n \sum_i \| E_i x \|_W : (E_i) \mbox{ is }S_n \mbox{
  allowable} \}\\
 & & \hskip .5in  \vee \sup_{n \geq s+1} 
  \sup \{ \theta_n \sum_i \| E_i x \|_W : (E_i) \mbox{ is }S_n \mbox{
  admissible} \},
\end{eqnarray*}
respectively. These norms  can also be defined as limits of appropriate
sequences. For $x \in c_{00}$ let
$$
\| x \|_{V,0}=\| x \|_{W,0}=\| x \|_{\infty}
$$
and for $m \in \N$ define:
\begin{eqnarray*}
\| x \|_{V,m+1} & = & \| x \|_\infty \vee \sup_n 
  \sup \{ \theta_n \sum_i \| E_i x \|_{V,m} : (E_i) \mbox{ is }S_n \mbox{
  allowable} \},\\
\| x \|_{W,m+1} & = & \| x \|_\infty \vee \sup_{n \leq s} 
  \sup \{ \theta_n \sum_i \| E_i x \|_{W,m} : (E_i) \mbox{ is }S_n \mbox{
  allowable} \}\\
 & & \hskip .5in  \vee \sup_{n \geq s+1} 
  \sup \{ \theta_n \sum_i \| E_i x \|_{W,m} : (E_i) \mbox{ is }S_n \mbox{
  admissible} \},
\end{eqnarray*}
Then 
$$
\| x \|_{V} =\lim_m\| x \|_{V, m},
\quad  \| x \|_{W} = \lim_m\| x \|_{W, m}. 
$$

\bigskip
\noindent
Then one can construct the  spaces 
 $V'= \Sigma( \alpha_n
\|. \|_{V,n})$, and $W'= \Sigma( \alpha_n \|. \|_{W,n})$. We show
 that $V'$ and $W'$ are
 $\ell_1$-saturated asymptotic  $\ell_1$ spaces for
vectors with disjoint supports.

\bigskip
\noindent
It is known \cite{CO} that if $\theta_n= \delta^n$ for some $\delta \in (0,1)$ then one 
can replace the ``allowable'' by ``admissible'' in the definition of $\| . \|_V$ to obtain
an equivalent norm for $V$. 
For this choice of $(\theta_n)$ the variant of the norm $\|.\|_{V,m+1}$ by
replacing
``allowable'' by ``admissible''  can be minorized and majorized up to a
uniform multiplicative constant by the norms $\|.\|_{V,m}$ and
$\|.\|_{V,m+1}$ respectively. This is enough to conclude that the new norms
lead to an equivalent norm for $V'$ (see   \cite{B}).

\section{The main Theorems}

\bigskip
\noindent
The following results is the main tool of our paper for constructing
weak Hilbert spaces.

\begin{Thm} \label{main}
If $X$ is an asymptotic $\ell_2$ space for vectors with disjoint
supports  then $X$ is a weak Hilbert space. 
\end{Thm}

\noindent
{\bf Proof}
Since the ideas needed for the proof of Theorem 3.1 exist in the
literature, we will just outline the proof.  Recall that the {\bf  
fast growing
hierarchy} is a sequence of functions on the natural numbers  
$(g_{n})$ which
is defined inductively by:  $g_{0}(n) = n+1$, and for $i\ge 0$,  
$g_{i+1}(n) =
g_{i}^{n}(n)$ where $g^{n}$ is the n-fold iteration of $g$ and  
$g^{0} = I$.

\bigskip
\noindent
{\bf Step I}
\begin{quote}
If $X$ is asymptotic-${\ell}_{2}$ with constant $C$ for vectors with disjoint supports, 
then for every  $i\ge 0$
any $g_{i}(n)$ normalized disjointly supported vectors with  
supports after n
 are $C^{i}$-equivalent to the unit vector basis of ${\ell}_{2}$.
\end{quote}
\noindent
We proceed by induction on $i$ with the case $i=0$ being trivial.  So, 
assume Step I holds for some $i\ge 0$ and let $\{x_{k}:1 \leq k \leq  
g_{i+1}(n)\}$
be a sequence of disjointly supported vectors in $X$ with supports  
after n.
For $1\le j \le n$ let
$$
E_{j} = \{k: g_{i+1}^{j-1}(n)\le k \le g_{i+1}^{j}-1\}.
$$
Then,
$$
\|\sum_{k=1}^{g_{i+1}^{n}(n)} x_{k}\| \stackrel{C}{\approx}
\left ( \sum_{j=1}^{n}\|\sum_{k\in E_{j}}x_{k}\|^{2}\right )^{1/2}
$$
Applying the induction hypotheses to each sum on the right we continue
this equivalence as
$$
\stackrel{C^{i+1}}{\approx} \left ( \sum_{k=1}^{g_{i+1}^{n}(n)}  
\|x_{k}\|^{2} \right )^{1/2}.
$$

\bigskip
\noindent
{\bf Step II}
\begin{quote}
If $X$ is asymptotic-${\ell}_{2}$ with constant $C$ for vectors
with disjoint supports then every
$n$-dimensional subspace of $X$ supported after n is 8$C^{3}$-isomorphic
to a Hilbert space and 8$C^{3}$-complemented in $X$.
\end{quote}
\noindent
If $E$ is a $5^{(5^{n})}$-dimensional subspace of $X$ supported  
after n, then
by a result of Johnson (See Proposition V.6 of \cite{CS}) there is
a subspace $G$ of $X$ spanned by $\le g_{3}(n)$ disjointly supported
vectors supported after n and an operator $V:E\rightarrow G$ with
$\|Vx-x\|\le \frac{1}{2}\|x\|$, for all $x\in E$.  Now, by Step I, we
have that $E$ is 2$C^{3}$-isomorphic to a Hilbert space.  It follows
\cite{J2} that every $5^{n}$-dimensional space of $X^{*}$ supported  
after n
is 4$C^{3}$-isomorphic to a Hilbert space and 4$C^{3}$-complemented
in $X^{*}$.  Therefore, every n-dimensional subspace of $X$ supported
after n is 8$C^{3}$-isomorphic to a Hilbert space and  
8$C^{3}$-complemented
in $X$.

\bigskip
\noindent
{\bf Step III}
\begin{quote}
Every asymptotic-${\ell}_{2}$ space for vectors with disjoint supports
is a weak Hilbert space.
\end{quote}
\noindent
If $E$ is a $2n$-dimensional subspace of $X$, let
$F =: E\cap (\spn_{k\ge n}e_{k})$.  Then $F$ is supported after n and
dim $F \ge n$ implies $F$ is $K$-isomorphic to a Hilbert space and
$K$-complemented in $X$ by Step II, where $K=8C^3$.  It follows from Definition \ref{weakH}
that $X$ is a weak Hilbert space. \hfill $\Box$

The $2$-convexification of certain Tsirelson spaces for obtaining weak
Hilbert spaces was first used in \cite{ADKM}. More generally we have
the following:

\begin{Cor} \label{ell_1,2}
If $X$ is an asymptotic $\ell_1$ space for vectors with disjoint
supports then $X^{(2)}$ is a weak Hilbert space. 
\end{Cor}

\noindent
{\bf Proof:} Let $(X,
\| . \|)$ be an asymptotic $\ell_1$ 
space for vectors with disjoint supports. Then
there exists $C>0$ such that 
for every  sequence of vectors $(x_i)$
with $(\supp x_i)$ being $S_1$-allowable, we have that 
$C \sum \| x_i \| \leq \| \sum x_i \|$. It suffices to prove that $X^{(2)}$ 
is an asymptotic $\ell_2$ space for vectors with disjoint support.
Let $(y_i)$ be a sequence of
vectors in $X^{(2)}$ with $(\supp y_i)$ being $S_1$-allowable. Then
$$
C^{1/2}( \sum \| y_i \|_{(2)}^2)^{1/2} =  C^{1/2} ( \sum \| y_i ^2
\|)^{1/2} 
 \leq  \| \sum y_i^2 \|^{1/2}
 =  \| ( \sum y_i^2)^{1/2} \|.
$$
Also, 
$$
\| (\sum y_i^2)^{1/2} \|_{(2)} =  \| \sum y_i^2 \|^{1/2}
 \leq  ( \sum \| y_i ^2 \|)^{1/2}= (\sum \| y_i \|_{(2)}^2)^{1/2}.
$$
\hfill $\Box$

\bigskip
\noindent
The spaces $V$, $W$, $V'$ and $W'$ are asymptotic $\ell_1$ spaces for
vectors with disjoint supports. Indeed, this is obvious for $V$ and
$W$. To see this for $V'$ let $n \in \N$ and vectors $(x_i)_{i=1}^n$ with disjoint
supports with $n \leq x_i$ for all $i$. Then:
\begin{eqnarray*}
\| \sum_{i=1}^n x_i \| & = & \sum_{m=1}^\infty \alpha_m 
\| \sum _{i=1}^n x_i \|_{V,m}
\geq  \sum_{m=1}^\infty \alpha_m \theta_1 \sum_{i=1}^n 
\| x_i \|_{V,m-1}\\
& \geq & \theta_1 u \sum_{i=1}^n 
\sum_{m=1}^\infty \alpha_{m-1} \| x_i \|_{V,m-1}
\geq  \theta_1 u \sum_{i=1}^n \| x_i \|_{V'}
\end{eqnarray*}The proof for $W'$ is
similar. Thus $V^{(2)}$, $W^{(2)}$, ${V'}^{(2)}$
and ${W'}^{(2)}$ are weak Hilbert spaces.

\begin{Prop} \label{asl1}
The spaces $V$ and $W$ do not contain an isomorph of $\ell_1$.
The spaces $V'$ and $W'$ are $\ell_1$-saturated without being
isomorphic to $\ell_1$. 
\end{Prop}

\noindent
The following Lemma will be used in the proof of Proposition
\ref{asl1}.

\begin{Lem} \label{c0sat}
For every $m \in \N$ the completion of $(c_{00}, \| . \|_{V,m})$ is a
$c_0$-saturated space.
\end{Lem}

\noindent
In order to prove this Lemma we need a result of Fonf along with
the notion of the boundary.

\begin{Def}
A subset $B$ of the unit sphere of the dual of a Banach space $X$ is
called a boundary for $X$ if for every $x \in X$ there exists $f \in
B$ such that $f(x)= \| x \|$.
\end{Def}

\begin{Thm}(\cite{F1}, see also \cite{F2}, \cite{H} \cite{DGZ})
Every Banach space with a countable boundary is $c_0$-saturated.
\end{Thm}

\noindent
{\bf Proof of Lemma \ref{c0sat}} Define inductively on $i \leq m$ the sets
$K^i$ of the unit ball of the dual of $(c_{00}, \| . \|_{V,m})$. Let
$K^0 = \{ \pm e_n: n \in \N \}$. For $i<m$ if $K^i$ has been defined then let
$$
K^{i+1}= K^i \cup \{ \theta_k (f_1 + \cdots + f_r): F_j \in K^m, \mbox{
for all }j, ( \supp f_j )_{j=1}^r \mbox{ is } S_k \mbox{ allowable } k=1,2,
\ldots \}.
$$
Then $K^m$ is a norming set for $(c_{00}, \| . \|_{V,m})$:
$$
\| x \|_{V,m}= \sup \{ |f(x)|: f \in K^m \}.
$$
It is easy to see that $K^m \cup \{ 0 \}$ is a
pointwise closed set since each $S_k$ is pointwise closed and $\lim_k
\theta_k =0$. The previous Theorem of Fonf finishes the proof of the
Lemma. \hfill $\Box$

\bigskip
\noindent
{\bf Proof of Proposition \ref{asl1}} 
The statement for $V$ and $W$ is obvious since they are reflexive \cite{AD}, 
\cite{ADKM}.\\
{\em $V'$ is $\ell_1$ saturated:} Let $(x_i)$ be an arbitrary block
basis of $V'$. It is enough to construct a normalized (in $V'$) block
basis $(v_i)$ of $(x_i)$ and an increasing sequence of positive
integers $1=p_1< p_2< p_3 < \cdots$ such that
$$
\sum_{m=p_i}^{p_{i+1}-1} \alpha_m \| v_i \|_{V,m} \geq \frac{1}{2}
$$ 
for all $i$. Once this is done then for $(\lambda_i) \in c_{00}$
\begin{eqnarray*}
\| \sum_{i=1}^n \lambda_i v_i \|_{V'} & = & \sum_{m=1}^\infty \alpha_m
 \| \sum_{i=1}^n \lambda_i v_i \|_{V,m} 
 =  \sum_{j=1}^\infty \sum_{m= p_j}^{p_{j+1}-1} \alpha_m 
\| \sum_{i=1}^n \lambda_i v_i \|_{V,m}\\
& \geq & \sum_{j=1}^\infty \sum_{m= p_j}^{p_{j+1}-1} \alpha_m 
\| \lambda_j v_j \|_{V,m}
 \geq  \frac{1}{2} \sum_{j=1}^\infty | \lambda_j | 
\end{eqnarray*} 
which shows that $(v_i)$ is equivalent to the unit vector basis of
$\ell_1$. In order to choose such $(v_i)$ and $(p_i)$ we use that for
every $m \in \N$ the norms $\| . \|_{V,m}$ and $\| . \|_{V,m+1}$ are
not equivalent. Thus for every $m,M,K \in \N$ there is $u$ in the span of
$(x_i)$ with 
$$
M \leq u, \quad \| u \|_{V,m} < \frac{1}{4}, \mbox{ and } 
\| u \|_{V,m+1} \geq K.
$$ 
Then 
$$
\sum_{i=1}^m \alpha_i \| u \|_{V,i} < \frac{1}{4} \mbox{ and }
\sum_{i=1}^\infty \alpha_i \| u \|_{V,i} \geq \alpha_{m+1}K.
$$
Let $v = \frac{u}{\| u \|_{V'}}$. By taking $K$ large enough we can
assume that 
$$
\sum_{i=1}^m \alpha_i \| v \|_{V,i} < \frac{1}{4}.
$$
Also choose $m' >m$ with
$$
\sum_{i=m'+1}^\infty \alpha_i \| v \|_{V,i} < \frac{1}{4}.
$$
Thus
$$
\sum_{i=m+1}^{m'} \alpha_i \| v \|_{V,i} \geq \frac{1}{2}.
$$ 
It only remains to show that for every $m \in \N$ the norms 
$\| . \|_{V,m}$ and $\| . \|_{V,m+1}$ are not equivalent on the span
of $(x_i')$. By the previous Lemma there is a block sequence $(y_i)$
of $(x_i)$ such that $\| y_i \|_{V,m}=1$ for all $i$ and 
$( (y_i), \| . \|_{V,m})$ is $2$-equivalent to the unit vector basis
of $c_0$. For $n \in \N$ let $k \in \N$ with 
$n \leq y_{k+1} < y_{k+2} < \cdots < y_{k+n}$. Thus 
$$
\| \sum_{i=k+1}^{k+n} y_i \|_{V,m} \leq 2
$$
yet
$$
\| \sum_{i=k+1}^{k+n} y_i \|_{V,m+1} \geq \theta_1 n.
$$
This proves the result.\\
{\em $W'$ is $\ell_1$-saturated:} Similar.\\
{\em $V'$ is not isomorphic to $\ell_1$:} If the statement were false
then the basis of $V'$ would be isomorphic to the unit vector basis of
$\ell_1$ (since every normalized unconditional basic sequence in $\ell_1$
is equivalent to the usual unit basis of $\ell_1$, \cite{LP}). 
Observe that for $x \in c_{00}$,
$$
\| x \|  =  \sum_m \alpha_m \| x \|_{V,m} 
 \leq  \sup_m  \| x \|_{V,m} = \| x \|_V.
$$
By \cite{ADKM} the  norm of $V$ can become arbitrarily smaller that the
$\ell_1$ norm on certain vectors. \\
{\em $W'$ is not isomorphic to $\ell_1$:} Similar. \hfill $\Box$

\begin{Rmk} \label{ell2saturated}
Note that a space $X$ with an unconditional basis contains $\ell_1$ if
and only if $X^{(2)}$ contains $\ell_2$. Since $V$ and $W$ do not contain an isomorph of
$\ell_1$, we obtain that $V^{(2)}$ and $W^{(2)}$ are weak Hilbert spaces 
which do not contain an isomorph of $\ell_2$. Since the spaces $V'$
and $W'$ are $\ell_1$ saturated without being isomorphic to $\ell_1$,
we obtain that the spaces ${V'}^{(2)}$, and
${W'}^{(2)}$ are $\ell_2$-saturated weak Hilbert spaces which are not
isomorphic to  $\ell_2$.
\end{Rmk}

\noindent
Thus the essential properties of the space $E_\alpha$ constructed by Edgington are
shared by ${V'}^{(2)}$ and ${W'}^{(2)}$.

\begin{Thm} \label{nonisomorphic}
Let $(\theta_n) \subset (0,1)$ with $\lim_n \theta_n^{1/n}=1$, $s \in
\N$,  and
$\beta=(\beta_n) \subset (0,1)$ with $\sum_n \beta_n =1$ and $0 < \inf
\frac{\beta_{n+1}}{\beta_n} \leq \sup
\frac{\beta_{n+1}}{\beta_n} <1$. Then ${V'}^{(2)}$ and
${W'}^{(2)}$ are not isomorphic to $E_\beta$.
\end{Thm}

\noindent
{\bf Proof:} Let $T:X \to E_{\beta}$ be
an isomorphism where $X$ is either ${V'}^{(2)}$ or
${W'}^{(2)}$. 
Since $T$ is an isomorphism  there exists $C>0$ such that
$$
\frac{1}{C} \|Tx \|_{E_\beta} \leq \| x \|_X \leq C \| Tx \|_{E_\beta}
$$
for all $x \in c_{00}$.  
Also, by \cite{E} (proof of Theorem 7) there exists $\delta >0$ such that
$$
\| Tx \|_{E_\beta} \leq C \| Tx \|_{T^{(2)}(\delta, S_1)}. 
$$
Thus for $x \in c_{00}$
\begin{equation} \label{norminequality}
\| x \|_X \leq C^2 \| Tx \|_{T^{(2)}(\delta, S_1)}.
\end{equation}
Since the unit vector basis $(e_i)$ of $X$ is weakly null, we can select a
subsequence $(e_{k_i})$ of $(e_i)$, a block sequence $(u_i)$ in $T^{(2)}(\delta, S_1)$
with non-negative coefficients and a number $K>0$ such that:
$$ 
\| T(e_{k_i}) - u_i \|_{T^{(2)}(\delta, S_1)} < \frac{\e}{2^i} \mbox{ and } 
\frac{1}{K} \leq \| u_i \|_{T^{(2)}(\delta, S_1)} \leq K \mbox{ for all }i,
$$
where $\e>0$ will be chosen later. Let $n \in \N$ to be selected
later. Let $(x_i)_{i \in I} \subset (0,1)$ for some $I \in S_n$ so that 
$$
\sum_{i \in I}x_i^2 =1,  \mbox{ and }
\| \sum_{i \in I} x_i^2 \frac{u_i^2}{\| u_i ^2 \|_{T(\delta,S_1)}} \|_{T(\delta,S_1)} 
\leq \delta^n +\e
$$
(\cite{OTW} Theorem 5.2 (a)). Then
\begin{eqnarray*}
\| T( \sum_{i \in I} x_i e_{k_i}) \|_{T^{(2)}(\delta, S_1)} & \leq & 
\| \sum_{i \in I} x_i u_i \|_{T^{(2)}(\delta, S_1)} + 
\sum_{i \in I} x_i \| Te_{k_i} - u_i \|_{T^{(2)}(\delta, S_1)}\\
&\leq & \| \sum_{i \in I} x_i^2 u_i^2 \| _{T(\delta, S_1)}^{1/2} + \e \\
&\leq & K \| \sum_{i \in I} x_i^2 \frac{u_i^2}{\| u_i ^2 \|_{T(\delta,S_1)}} 
\| _{T(\delta, S_1)}^{1/2} + \e \\
& \leq & K (\delta^n+ \e)^{1/2} + \e. 
\end{eqnarray*}
On the other hand if $Y=V'$ when $X={V'}^{(2)}$ or $Y=W'$ when
$X={W'}^{(2)}$ then 
$$
\| \sum_{i \in I} x_i e_{k_i} \|_X  =  
 \| \sum_{i \in I} x_i^2 e_{k_i} \|_Y ^{1/2} \geq 
( \sum_{m=1}^\infty \beta_m  \theta_n \sum_{i \in I} x_i^2 )^{1/2}
= \sqrt{ \theta_n}
$$
Therefore (\ref{norminequality}) gives
$$
\sqrt{\theta_n} \leq C^2(K (\delta^n + \e)^{1/2}+ \e).
$$
But since $\lim_n \theta_n^{1/n}=1$, $n$ and $\e$ can be chosen
so that this  inequality fails.
\hfill $\Box$

\end{document}